\documentclass[reqno,12pt]{amsart}

\usepackage{epsf}
\usepackage{amssymb}
\usepackage{amsmath}
\usepackage{longtable}

\date{}

\def\pic#1#2{ \def\epsfsize##1##2{0.7##1}\raisebox{#1}{$\vcenter{\hbox{\epsffile{#2.eps}}}$}}

\def\tablepic#1#2{ \def\epsfsize##1##2{0.5##1}\raisebox{#1}{$\vcenter{\hbox{\epsffile{#2.eps}}}$}}

\def\minipic#1#2{ \def\epsfsize##1##2{0.3##1}\raisebox{#1}{$\vcenter{\hbox{\epsffile{#2.eps}}}$}}

\theoremstyle{plain}
\newtheorem{theorem}{Theorem}
\newtheorem*{sbi}{Slice-Bennequin Inequality}
\newtheorem*{mo}{Theorem (Morton \cite{Mo})}
\newtheorem*{corollary}{Corollary}
\newtheorem*{lemma}{Lemma}
\newtheorem*{proposition}{Proposition}

\theoremstyle{remark}
\newtheorem*{definition}{\scshape{Definition}}

\newtheorem*{examples}{\scshape{Examples}}
\newtheorem*{remark}{\scshape{Remark}}
\newtheorem*{remarks}{\scshape{Remarks}}
\newtheorem*{ack}{\scshape{Acknowledgements}}

\def\Z{{\mathbb Z}}

\def\R{{\mathbb R}}

\author{Sebastian Baader}
\title{Slice and Gordian numbers of track knots}

\begin{document}

\begin{abstract}
We present a class of knots associated with labelled generic immersions of
intervals into the plane and compute their Gordian numbers and 4-dimensional
invariants. At least $10\%$ of the knots in Rolfsen's table belong to this
class of knots. We call them track knots. They are contained in the class of
quasipositive knots. In this connection, we classify quasipositive knots and
strongly quasipositive knots up to 10 crossings.
\end{abstract}

\maketitle

\section{Introduction}

Several classes of knots are closely related to generic immersions of compact
1-manifolds into the plane. The class of track knots we shall present
subsequently is a partial generalization of the class of divide knots.
A divide is the intersection of a plane curve with the unit disk in $\R^2$,
provided the plane curve is transverse to the unit circle. 
The concept of knots associated with divides is due to Norbert A'Campo and
emerged from the study of isolated singularities of complex plane curves
(see \cite{AC1}). In \cite{AC2} and \cite{AC3}, A'Campo specified some
properties of divide knots, including fiberedness and a Gordian number result.
Mikami Hirasawa gave an algorithm for drawing diagrams of divide links and
extended the Gordian number result to certain arborescent links (see
\cite{H}). A large extension of the class of divide links was introduced by
William Gibson and Masaharu Ishikawa \cite{GI}. They kept the Gordian number
result, too. Tomomi Kawamura \cite{Ka2} and Ishikawa independently proved the
quasipositivity of these links of free divides.

Borrowing from all these, we propose a new construction of knots associated
with labelled generic immersions of intervals into the plane. 

\smallskip
Let $C$ be the image of a generic immersion of the interval $[0,1]$ into the
plane. In particular, $C$ has no multiple points apart from a finite number of
transversal double points, none of which is the image of 0 or 1. Further we
enrich $C$, as follows (see Figure 1 for an illustration).

\begin{enumerate}
\renewcommand{\labelenumi}{(\roman{enumi})}
\item A small disk around each double point of $C$ is cut into four regions by
$C$. Label each of these regions by a sign, such that the sum of the four
signs is non-negative. There are four types of patterns of signs around a
double point, called $a$, $b$, $c$ and $d$. They are shown in Figure 2. If the
tangent space $T_pC$ at a double point $p$ of $C$ is the set $\{(x,y) \in \R^2
\mid (y-y(p))^2=(x-x(p))^2\}$, then we may represent patterns of four signs at
$p$ by one of the following symbols:

$a$, $a_1$, $b$, $b_1$, $b_2$, $b_3$, $c$, $c_1$, $c_2$, $c_3$, $d$.

\smallskip
An index $i$ at a symbol means that the corresponding pattern has to be turned
counter-clockwise by the angle $i \frac{\pi}{2}$. 

\smallskip
For example, $b_1$ \!\!\!\!\! \minipic{-0pt}{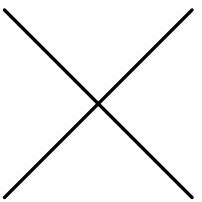} stands for the pattern
\minipic{-0pt}{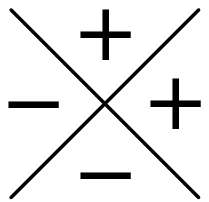} . 

\smallskip
Henceforth we shall use these symbols. 

\smallskip
\item Specify a finite number of different points $p_1$, $p_2$, \ldots, $p_r$
on the edges of $C$ (i.e. on the connected components of $C-$\{double
points\}, such that $C-\{p_1, \, p_2, \ldots, p_r \}$ is simply connected, but
not necessarily connected. $r$ is greater than or equal to the number of double points of $C$.
\end{enumerate}

A labelled generically immersed interval in the plane will always be denoted by $C_{\lambda}$.

\smallskip
\begin{figure}[ht]
\pic{-0pt}{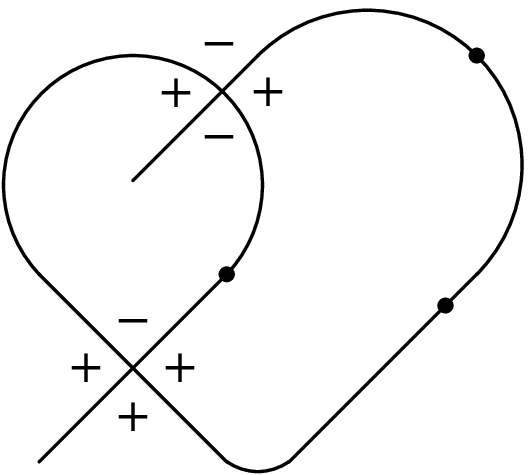}
\caption{}
\end{figure}

\smallskip
\begin{figure}[ht]
a \pic{-0pt}{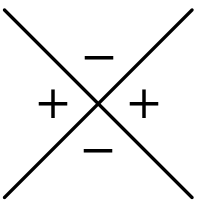} \qquad \qquad b \pic{-0pt}{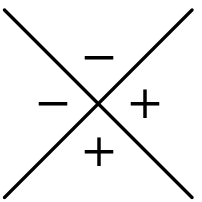}

\bigskip
c \pic{-0pt}{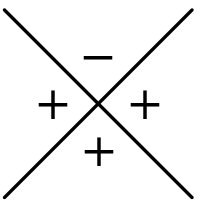} \qquad \qquad d \pic{-0pt}{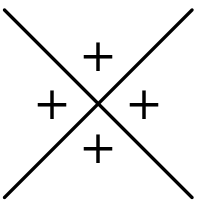}
\caption{Patterns of signs}
\end{figure}

The following algorithm associates a knot diagram, hence a knot in the
3-space, to a labelled generically immersed interval $C_{\lambda}$.

\begin{enumerate}
\item Draw a parallel companion of $C_{\lambda}$. In other words, replace
$C_{\lambda}$ by the boundary of a small band following $C_{\lambda}$. Join
the two strands with an arc at both end points of $C_{\lambda}$ and orient the
resulting plane curve clockwise, in regard of the small band (see Figure~3).

\begin{figure}[ht]
\pic{-0pt}{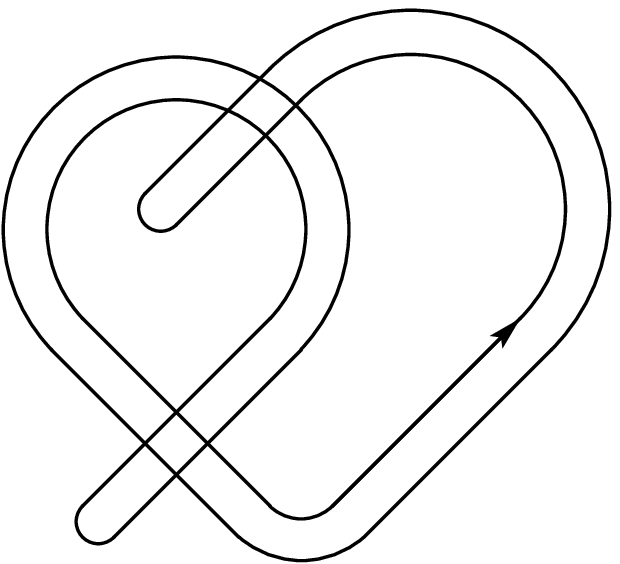}
\caption{}
\end{figure}

\item At each double point of $C_{\lambda}$, place over- and under-crossings
according to the signs of the four regions, as shown in Figure 4.

\smallskip
\begin{figure}[ht]
a \pic{-0pt}{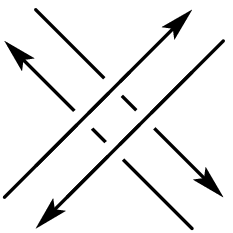} \qquad \qquad b \pic{-0pt}{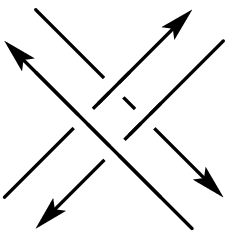}

\medskip
c \pic{-0pt}{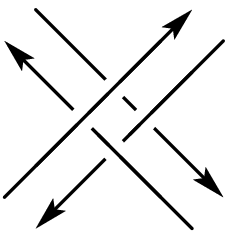} \qquad \qquad d \pic{-0pt}{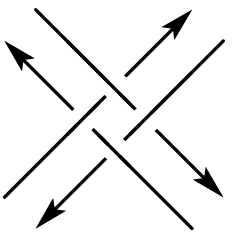}

\bigskip
\pic{-0pt}{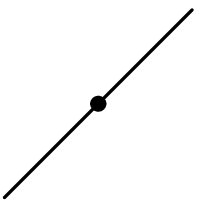} \quad $\longrightarrow$ \quad \pic{-0pt}{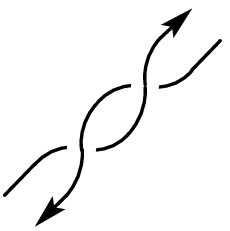}
\caption{}
\end{figure}

The characters a, b, c and d stand for 'above', 'between', 'conventional'
and 'double', respectively. 'Conventional' crossings appear in the visualization of links of divides, see Hirasawa \cite{H}.

\item Add a full twist to the band at each specified point of $C_{\lambda}$, in a manner that gives rise to two positive crossings (see Figure 4).
\end{enumerate}

The knot diagram arising from $C_{\lambda}$ by these three steps will be denoted by $D(C_{\lambda})$, the corresponding knot by $K(C_{\lambda})$.

\begin{definition} A \emph{track knot} is a knot which can be realized as a knot
associated with a labelled generically immersed interval $C_{\lambda}$. If it
can be realized without any double point of type $b$, then we call it a
\emph{special track knot}.
\end{definition}

\begin{remark}
We observe that the classes of track knots and special track knots are closed
under connected sum. The connected sum operation corresponds to the gluing of
two labelled immersed intervals along end points. This is not true for knots
of free divides; the connected sum of the free divide knot $5_2$ (in Rolfsen's
numbering \cite{Ro}) with itself is not a free divide knot.
\end{remark}

\section{Slice and Gordian Numbers}

Let $L$ be an oriented link with $n$ components in $S^3=\partial B^4$. The
slice number $\chi_s(L)$ of $L$ is the maximal Euler characteristic of all
smooth, oriented surfaces in $B^4$ which are bounded by $L$ and have no closed
components. The surfaces in consideration need not be connected. If $K$ is a
knot, the 4-genus $g^*(K)$ is defined as $\frac{1}{2}(1-\chi_s(K))$. The clasp
number $c_s(L)$ of a link $L$ is the minimal number of transversal double
points of $n$ generically immersed disks in $B^4$ with boundary $L$. We will
also be concerned with the Gordian number $u(L)$, which is the minimal number 
of crossing changes needed to transform $L$ into the trivial link with $n$
components. The following two inequalities relate these numbers:

\begin{equation}
u(L) \geqslant c_s(L) \geqslant \frac{1}{2}(1-\chi_s(L)).
\label{1}
\end{equation}

\noindent
They can be shown by purely geometrical arguments, see Kawamura \cite{Ka1}.

Gordian numbers and 4-dimensional invariants of track knots are easy to
determine. Let $K$ be a track knot associated with a labelled generically
immersed interval $C_{\lambda}$. Further let $A$, $B$, $C$ and $D$ be the
numbers of double points of $C_{\lambda}$ with patterns of signs of type $a$,
$b$, $c$ and $d$, respectively.

\begin{theorem} The clasp number and the 4-genus of $K$ equal $C+2D$. If $B$ is zero, then the Gordian number and the ordinary genus of $K$ equal $C+2D$, too.
\end{theorem}

\begin{corollary} Both the clasp number and the 4-genus are additive under
connected sum of track knots. Moreover, the Gordian number is additive under
connected sum of special track knots.
\end{corollary}

\begin{remark} The connected sum of a knot with its mirror image always
bounds an embedded disk in $B^4$, thus the clasp number and the 4-genus are
not additive under connected sum of knots in general. It is still a conjecture
that the Gordian number is additive under connected sum of knots (see M.~Boileau and C.~Weber \cite{BW}).
\end{remark}

\begin{proof}[\rm{Proof of Theorem 1}] We first show that the 4-genus of $K$
does not exceed $C+2D$. If $C=D=0$, then $K$ is clearly slice, i.e. $K$ bounds
a disk in $B^4$. Indeed, the band following $C_{\lambda}$ provides an immersed
disk in $S^3$ with boundary $K$. At each double point of type $b$ we may push
a part of the band into $\mathring{B}^4$ to get an embedded disk. But then, at
each double point of type $c$, we add one handle to the band, as Figure 5(c)
suggests. Similarly, we add two handles to the band at each double point of type $d$, see Figure 5(d). This creates an embedded surface in $B^4$ of genus $C+2D$ with boundary $K$. If $B=0$, it is an embedded surface in $S^3$.

\begin{figure}[ht]
c \pic{-0pt}{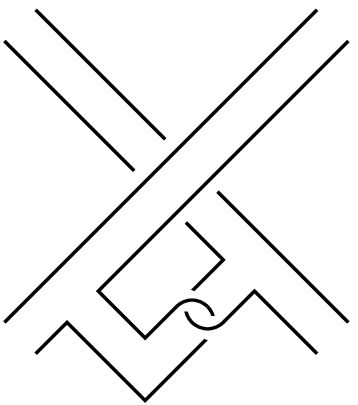} \qquad \qquad d \pic{-0pt}{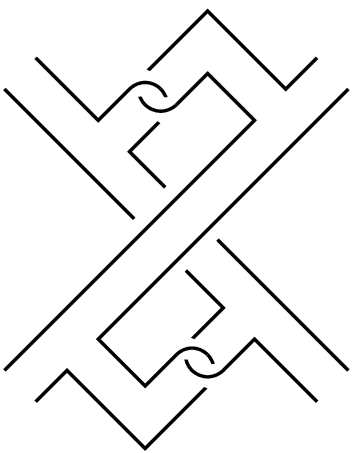}
\caption{}
\end{figure}

We remark that the spots where we add handles to the band can be interpreted
as clasp singularities of the immersed band. Therefore the clasp number of $K$
does not exceed $C+2D$, either. Next, we show that the Gordian number of $K$
does not exceed $C+2D$, provided $B$ is zero. If $C=D=0$, then $K$ is the
unknot since it bounds an embedded disk in $S^3$. On a knot diagram level,
double points of type $c$ differ from double points of type $a$ only by one
crossing change, see Figure 4. Similarly, double points of type $d$ differ
from double points of type $a$ by two crossing changes. Hence we conclude
$u(K) \leqslant C+2D$. 

\smallskip
We still have to prove that $C+2D$ is a lower bound for the four numbers in
question. If we prove $g_*(K) \geqslant C+2D$, then we are done, thanks to
(1). For this purpose we need the slice-Bennequin inequality. Let $D_L$ be the
diagram of an oriented link $L$. The writhe $w(D_L)$ is the number of positive
minus the number of negative crossings of the diagram $D_L$. Smoothing $D_L$
at all crossings produces a union of Seifert circles. Let $s(D_L)$ be their
number.

\begin{sbi} $\chi_s(L) \leqslant s(D_L)-w(D_L)$.
\end{sbi}

The slice-Bennequin inequality was first established for closed braid diagrams
by Lee Rudolph \cite{Ru1}; the proof of the general case can be found in
Rudolph \cite{Ru2} and Kawamura \cite{Ka1}. A '3-dimensional' version of the
inequality (concerning Seifert surfaces) was proved by Daniel Bennequin
\cite{Be}.

\smallskip
Now let us compute $w(D(C_{\lambda}))$ and $s(D(C_{\lambda}))$ for the knot
diagram of the labelled generically immersed interval $C_{\lambda}$.

\begin{enumerate}
\item $w(D(C_{\lambda}))=2C+4D+2r$, where $r$ is the number of specified points on $C_{\lambda}$.

\smallskip
\item Each double point and each specified point of $C_{\lambda}$ gives rise
to a small Seifert circle, see Figure 6. Moreover, each connected component of
$C_{\lambda}-\{p_1, \, p_2, \ldots, p_r \}$ gives one Seifert circle. The
number of connected components of $C_{\lambda}-\{p_1, \, p_2, \ldots, p_r \}$
being $1+r-(A+B+C+D)$, we conclude 

\medskip
\noindent
$s(D(C_{\lambda}))=A+B+C+D+r+1+r-(A+B+C+D)=2r+1$.

\end{enumerate}

Thus the slice-Bennequin inequality yields $\chi_s(K) \leqslant 1-2C-4D$ and
$g^*(K)=\frac{1}{2}(1-\chi_s(K)) \geqslant C+2D$.
\end{proof}

\begin{figure}[ht]
\pic{-0pt}{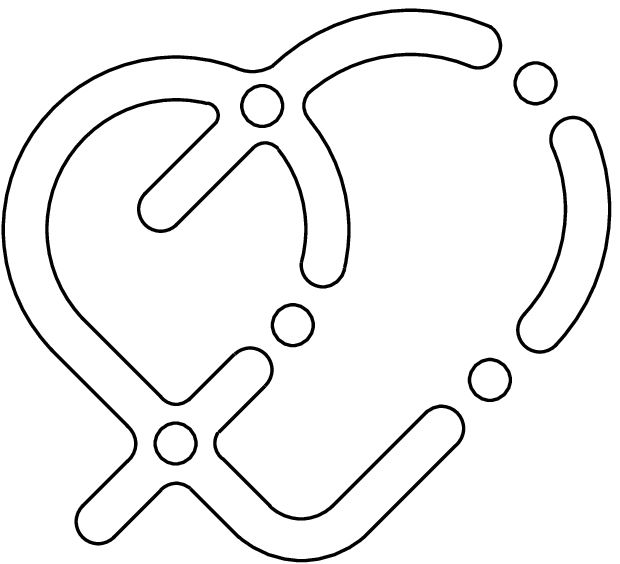}
\caption{}
\end{figure}

\vbox{%
\begin{remarks}\quad
\begin{enumerate}
\renewcommand{\labelenumi}{(\roman{enumi})}
\item If we renounce twisting the band at some specified points, then the
statements of Theorem 1 are no longer true. The labelled immersed interval
(without specified points) of Figure 7 has one double point of type $c$ and
gives the unknot.

\item The statement of Theorem~1 about the Gordian number can be extended for
track knots with $B=1$. However, if $B \geqslant 2$, then the Gordian number may be greater than $C+2D$. E.g. the knots $9_{46}$ and $10_{140}$ are slice track knots (see Table 2) and their Gordian numbers are certainly not zero.
\end{enumerate}
\end{remarks}
}

\begin{figure}[ht]
\pic{-0pt}{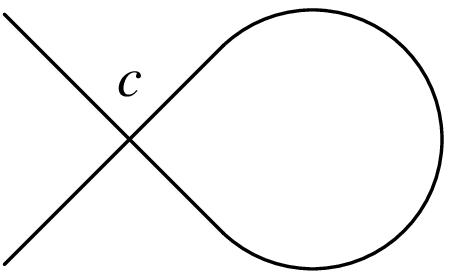}
\caption{}
\end{figure}

\section{The Knots $10_{131}$ and $10_{148}$}

The knot $10_{131}$ is the track knot corresponding to the labelled immersed
interval of Figure 8. Its 4-genus and Gordian number equal 1. The latter is
declared unknown in Akio Kawauchi's table of knots \cite{Kaw}. It is a curious
fact that we can see the unknotting operation on its minimal diagram both in
Rolfsen's and Kawauchi's table. This was already observed by Alexander
Stoimenow in \cite{St2}.

\begin{figure}[ht]
\pic{-0pt}{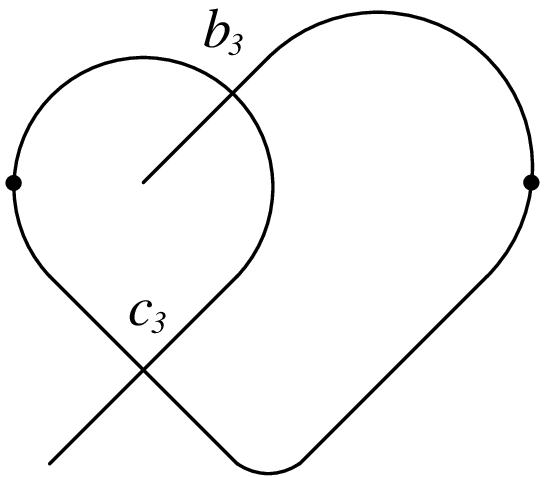}
\caption{}
\end{figure}

The knot $10_{148}$ is not a genuine track knot; it corresponds to a labelled
immersed interval with too little specified points, see Figure 9. Copying the
first part of the proof of Theorem 1, we see that its 4-genus is 1 at most.
Since the knot $10_{148}$ is already known not to be slice, we conclude that
its 4-genus is 1. This entry in Kawauchi's table of knots has been
corrected a few years ago, see \cite{Kaw2}. However, we cannot decide whether
its Gordian number is 1 or 2. 

\smallskip
\begin{figure}[ht]
\pic{-0pt}{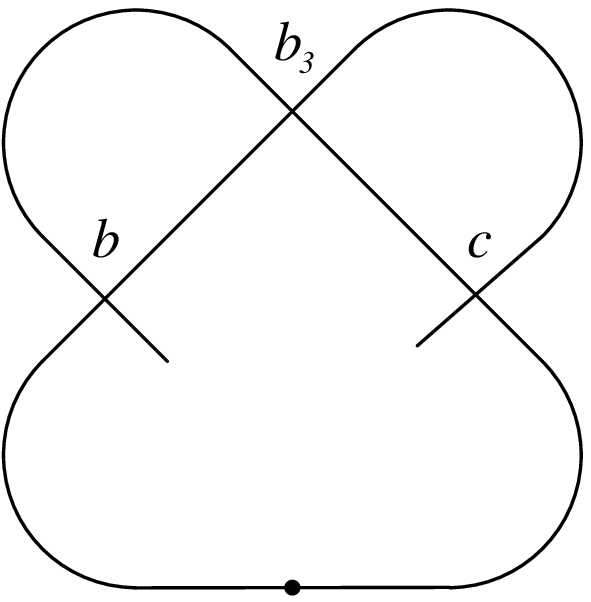}
\caption{}
\end{figure}

\section{The HOMFLY Polynomial and Quasipositivity}

The HOMFLY polynomial $P_L(v,z) \in \Z[v^{\pm 1},z^{\pm 1}]$ of an oriented link $L$ is defined by the following two requirements (see \cite{FYHLMO}):

\medskip
\noindent 
1. Normalization: $P_{O}(v,z)=1$,

where $O$ stands for the regular diagram consisting of one trivial circle.

\medskip
\noindent
2. Relation: $\frac{1}{v}P_{D_+}(v,z)-vP_{D_-}(v,z)=zP_{D_o}(v,z)$.

\smallskip
Here $D_+$, $D_-$ and $D_o$ denote regular diagrams which coincide
outside a standard disk and differ, as in Figure 10, inside this disk.

\begin{figure}[ht]
$D_+$ \pic{-0pt}{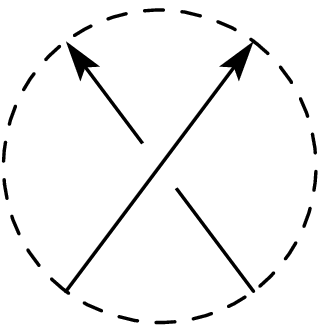} \quad \quad $D_-$ \pic{-0pt}{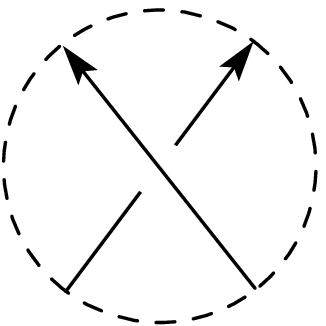} \quad \quad
$D_o$ \pic{-0pt}{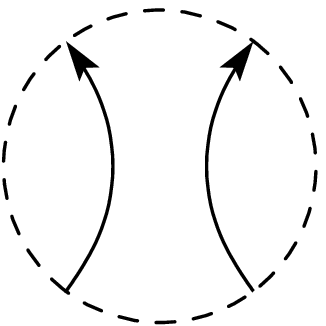}
\caption{}
\end{figure}

\medskip
Writing $P_L(v,z)=\sum_{k=e(L)}^{E(L)}a_k(z)v^k$, with $a_{e(L)}(z)$,
$a_{E(L)}(z) \neq 0$, as a Laurent polynomial in one variable $v$, we define
its range in $v$ as $[e(L),E(L)]$. H. R. Morton gave some bounds for $e(L)$
and $E(L)$ in terms of the writhe and the number of Seifert circles of a
diagram of $L$.

\begin{mo} For any diagram $D_L$ of an oriented link $L$
$$w(D_L)-(s(D_L)-1) \leqslant e(L) \leqslant E(L) \leqslant w(D_L)+(s(D_L)-1).$$
\end{mo}

The first inequality is tailor-made for track knots.

\begin{theorem} $2g^*(K) \leqslant e(K)$ for any track knot $K$.
\end{theorem}

\begin{proof}[\rm{Proof}] Choose a track knot diagram $D$ of $K$. The proof of
Theorem 1 tells us that $g^*(K)=\frac{1}{2}(1-s(D)+w(D))$, which is exactly half the lower bound in Morton's theorem.

\end{proof}

Theorem 2 draws our attention to quasipositive knots. A quasipositive knot is
a knot which can be realized as the closure of a quasipositive braid. A
quasipositive braid is a product of conjugates of a positive standard generator of the braid group. The slice-Bennequin inequality being an equality for closed quasipositive braid diagrams, we see that Theorem 2 is true both for track knots and for quasipositive knots. 

\begin{theorem} Track knots are quasipositive.
\end{theorem}

We adopt the pattern of Takuji Nakamura's proof of strong
quasipositivity of positive links (see \cite{N}). Any planar knot diagram
gives rise to a system of Seifert circles with signed arcs, where each arc
stands for a crossing joining two Seifert circles, as shown in Figure 11. The
sign of an arc tells us whether the crossing is positive or negative.

\begin{figure}[ht]
\pic{-0pt}{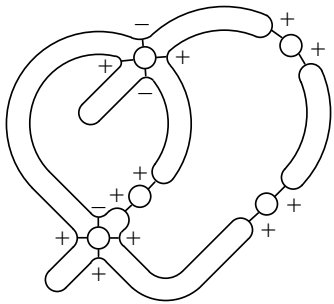}
\caption{A system of Seifert circles}
\end{figure}

\begin{definition} A knot diagram is \emph{quasipositive} if its set of
crossings can be partitioned into single crossings and pairs of crossings, such that the following three conditions are satisfied.

\begin{enumerate}
\item Each single crossing is positive.

\item Each pair of crossings consists of one positive and one negative
crossing joining the same two Seifert circles.

\item A pair of crossings does not separate other pairs of crossings. More
precisely, going from one crossing of a pair to its opposite counterpart along
a Seifert circle, one cannot meet only one crossing of a pair.
\end{enumerate}
\end{definition}

\vbox{%
\begin{examples}\quad
\begin{enumerate}
\item[$\bullet$] Positive knot diagrams are obviously quasipositive.

\item[$\bullet$] Track knot diagrams are quasipositive: negative arcs are
incident with a small Seifert circle corresponding to a double point of type
$a$, $b$ or $c$. They can be paired with neighbouring positive crossings of
the same small Seifert circle (see Figure 12). At this point, it is essential
that $C_{\lambda}-\{p_1, \, p_2, \ldots, p_r \}$ is simply connected. This
guarantees that pairs of crossings do not get entangled (see Figure 11). 

\item[$\bullet$] Quasipositive braid diagrams are quasipositive.
\end{enumerate}
\end{examples}
}

\begin{figure}[ht]
\pic{-0pt}{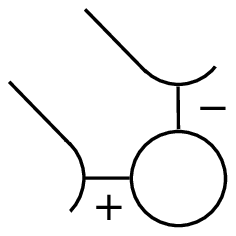}
\caption{A pair of crossings}
\end{figure}

\begin{lemma} A quasipositive knot diagram represents a quasipositive knot.
\end{lemma}

\begin{proof}[\rm{Proof}] Any link diagram can be deformed into a braid
representation, i.e. a system of concentric Seifert circles,  by a finite
sequence of bunching operations or concentric deformations of two types,
without changing the writhe and the number of Seifert circles of the link
diagram. This algorithm is due to Shuji Yamada, see \cite{Y}. We shall explain
these two deformations and their effect on quasipositive knot diagrams.

First of all, we may consider only knot diagrams which have an outermost
Seifert circle $S_1$, i.e. one that contains all the other Seifert circles.
This corresponds to choosing a point on the sphere $S^2$ appropriately.

If $S_1$ contains a maximal Seifert circle $S_2$ with the opposite orientation
of $S_1$, then we apply a concentric deformation of type I to $S_2$, as shown
in Figure~13.

\begin{figure}[ht]
\tablepic{-0pt}{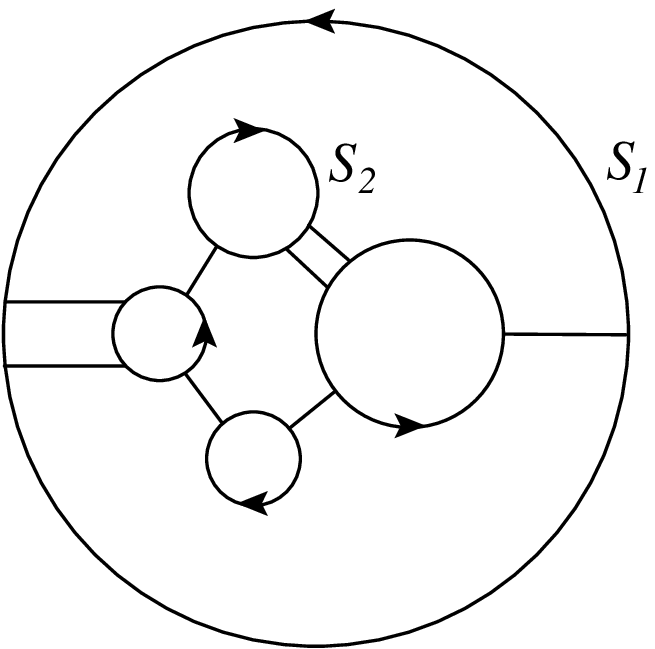} $\overset{c_I}\longrightarrow$
\tablepic{-0pt}{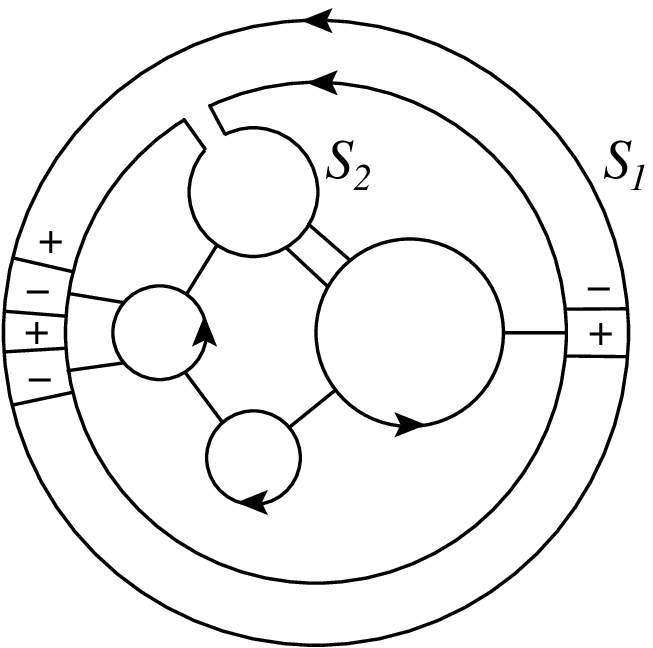}
\caption{A concentric deformation of type I}
\end{figure}

If $S_1$ contains maximal Seifert circles with the same orientation as $S_1$
only, then we apply a concentric deformation of type II to any of these
maximal Seifert circles, say to $S_2$, as shown in Figure 14.

\begin{figure}[ht]
\tablepic{-0pt}{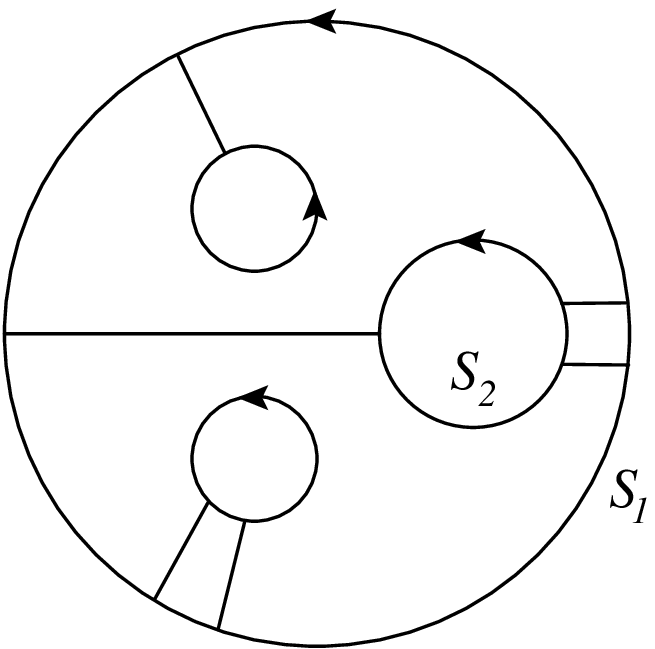} $\overset{c_{II}}\longrightarrow$
\tablepic{-0pt}{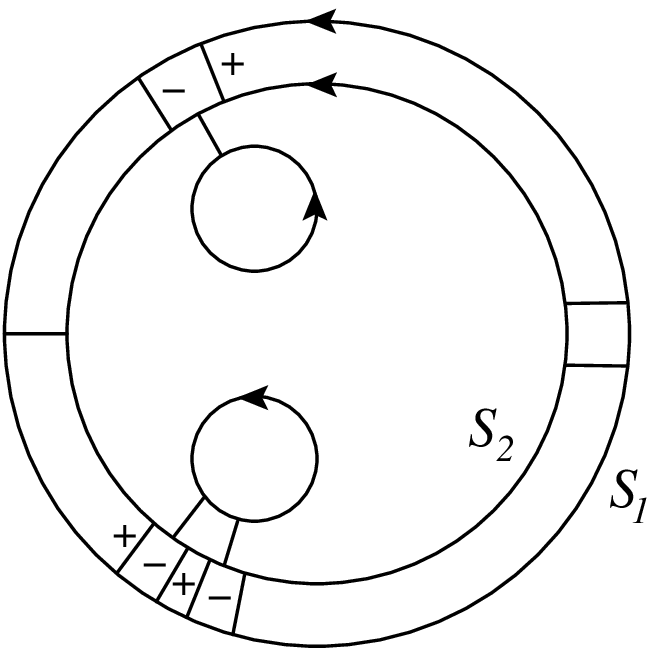}
\caption{A concentric deformation of type II}
\end{figure}

In the next step, we consider maximal Seifert circles inside $S_2$, and so on.
This algorithm clearly ends in a braid representation. Now we observe that
concentric deformations of both types preserve the quasipositivity of knot
diagrams in the above sense. They merely introduce new pairs of crossings,
which do not get entangled. Figures 15 and 16 show how a positive crossing (or
a pair of crossings, respectively) gets more 'conjugated' by new pairs of
crossings after a concentric deformation.

\begin{figure}[ht]
\tablepic{-0pt}{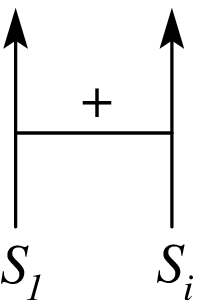} $\overset{c_*}\longrightarrow$
\tablepic{-0pt}{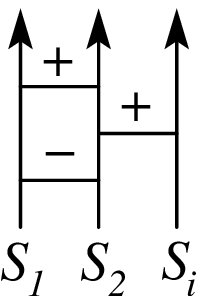}
\caption{}
\end{figure}

\begin{figure}[ht]
\tablepic{-0pt}{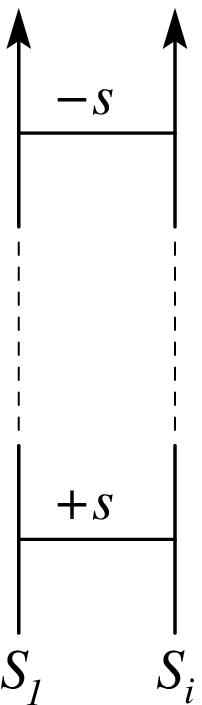} $\overset{c_*}\longrightarrow$
\tablepic{-0pt}{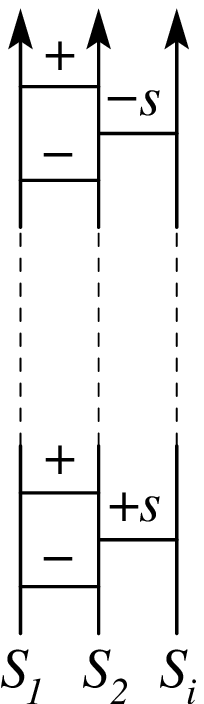}
\caption{}
\end{figure}

Thus, starting with a quasipositive knot diagram, we end up with a
quasipositive braid diagram, which clearly represents a quasipositive knot.

\end{proof}

Theorems 2 and 3 reduce the number of potential track knots. In the following,
we consider prime knots up to 10 crossings. Looking at Kawauchi's table of
knots, we see that 60 of 249 prime knots up to 10 crossings satisfy the
inequality $2g^*(K) \leqslant e(K)$. Among these 60 knots, 42 have positive
diagrams:

\medskip
\noindent $3_1$, $5_1$, $5_2$, $7_1$, $7_2$, $7_3$, $7_4$, $7_5$, $8_{15}$,
$8_{19}$, $9_1$, $9_2$, $9_3$, $9_4$, $9_5$, $9_6$, $9_7$, $9_9$, $9_{10}$,
$9_{13}$, $9_{16}$, $9_{18}$, $9_{23}$, $9_{35}$, $9_{38}$, $9_{49}$,
$10_{49}$, $10_{53}$, $10_{55}$, $10_{63}$, $10_{66}$, $10_{80}$, $10_{101}$,
$10_{120}$, $10_{124}$, $10_{128}$, $10_{134}$, $10_{139}$, $10_{142}$,
$10_{152}$, $10_{154}$, $10_{161}$.

\vbox{%
\begin{remarks}\quad
\begin{enumerate}
\renewcommand{\labelenumi}{(\roman{enumi})}
\item Since knots are always listed up to mirror image, we must be more precise: 'a knot $K$ satisfies the inequality $\ldots$' means 'either $K$ or its mirror image $!K$ satisfies the inequality $\ldots$'.

\item The 4-genus of the knot $10_{51}$ is not known. However, it is known not
to be slice, hence the inequality $2g^*(10_{51}) \leqslant e(10_{51})=0$ is not
satisfied.
\end{enumerate}
\end{remarks}
}

According to the result of Nakamura \cite{N} and Rudolph \cite{Ru2}, positive
knots are strongly quasipositive, i.e. they can be realized as the closure of a
braid which is a product of positive embedded bands of the form

$$\sigma_{i,j}=(\sigma_i \cdots \sigma_{j-2})\sigma_{j-1}(\sigma_i \cdots
\sigma_{j-2})^{-1},$$

where $\sigma_i$ is the i-th positive standard generator of the braid group.
The remaining 18 knots are listed in Table~1, except for the knot $10_{132}$,
which is not quasipositive. Stoimenow already pointed out that the
quasipositivity of the knot $10_{132}$ would imply the quasipositivity of its
untwisted 2-cable link, together with a violation of Morton's inequality,
which is a contradiction (see \cite{St1}). Table~1 contains one strongly
quasipositive, non-positive knot: $10_{145}$. It is non-positive since it is
non-homogeneous (see P.~R.~Cromwell \cite{C}). The other 16 knots are not
strongly quasipositive since their 4-genus is smaller than their genus. In
particular, they are non-positive. So in Table~1 we list all quasipositive,
non-positive prime knots up to 10 crossings in Rolfsen's numbering, together
with a quasipositive braid representation, the 4-genus $g^*$ and the ordinary
genus $g$. In the second column a, b, $\ldots$ and A, B, $\ldots$ stand for
$\sigma_1$, $\sigma_2$, $\ldots$ and $\sigma_1^{-1}$,
$\sigma_2^{-1}$, $\ldots$ and have nothing to do with symbols of labelled
immersed intervals. Parentheses should help to recognize positive bands. The
braid of the knot $10_{145}$ is strongly quasipositive.

This classification of quasipositive and strongly quasipositive knots gives us
an interesting criterion for detecting strongly quasipositive knots.

\begin{proposition} A knot with 10 crossings at most is strongly quasipositive, if and only if it is quasipositive and its 4-genus equals its ordinary genus.
\end{proposition}
 
\pagebreak
\begin{longtable}{| c | l | c | c |}
\caption{Quasipositive, non-positive prime knots up to 10 crossings} \\
\hline
\endhead
\hline
Knot  &  quasipositive braid representation  &  $g^*$  &  $g$  \\ \hline 
$8_{20}$  & (abAbaBA)(baB)  & 0 & 1  \\ \hline
$8_{21}$  & (abA)b(Abba)  & 1 & 2  \\ \hline
$9_{45}$  & a(Bcb)b(bacB)  & 1 & 2  \\ \hline
$9_{46}$  & (abbcBBA)(bacB)  & 0 & 1  \\ \hline
$10_{126}$  & aa(aaabAAA)b  & 1 & 3  \\ \hline
$10_{127}$  & abbb(bAbbaB)  & 2 & 3  \\ \hline
$10_{131}$  & a(aaBCbdBcbAA)(BcbdcBCb)d(Bcb)  & 1 & 2  \\ \hline
$10_{133}$  & aab(bDCbcdB)(bCBcACbcdCBcaCbcB)(bCBcaCbcB) & 1 & 2  \\ \hline
$10_{140}$  & (abbbcBBBA)b(Cbc)  & 0 & 2  \\ \hline
$10_{143}$  & a(BBBaaabbb)  & 1 & 3  \\ \hline
$10_{145}$  & (abA)cd(abA)(bcB)(bcdCB)(cdC)b  & 2 & 2  \\ \hline
$10_{148}$  & ab(bbacBB)(cbC)  & 1 & 3  \\ \hline
$10_{149}$  & a(bbCbccBB)a(bcccB)  & 2 & 3  \\ \hline
$10_{155}$  & (abA)(ABcbCba)(bcB)  & 0 & 3  \\ \hline
$10_{157}$  & a(Baab)b(baaB)  & 2 & 3  \\ \hline
$10_{159}$  & a(BBaabb)(baB)  & 1 & 3  \\ \hline
$10_{166}$  & (abcBA)(acbA)(Bcb)(Aba)  & 1 & 2  \\ \hline
\end{longtable}

We conclude this section with some questions and problems arising from the
study of track knots and quasipositive knots.

\begin{enumerate}

\item Does there exist a quasipositive knot which is not a track knot?

\noindent (The free divide knots $9_{16}$, $10_{124}$, $10_{152}$ and $10_{154}$ might be good candidates.)

\smallskip
\item Classify track knots up to 10 crossings. For this purpose, find new
criterions for detecting track knots.
\smallskip
\item Is it true that a knot is strongly quasipositive, if and only if it is
quasipositive and its 4-genus equals its ordinary genus? In particular, is it
true that special track knots are strongly quasipositive?

\smallskip
\item Do alternating quasipositive knots have positive diagrams? (Up to ten
crossings, this is true.)

\smallskip
\item Generalize the class of track knots in order to get some new Gordian
number results.

\smallskip
\item Prove the additivity of the Gordian number under connected sum of knots.

\end{enumerate}

\section{Examples of Track Knots}

In this section we look at the labelled immersed interval shown in Figure 17.
It has two double points and two specified points.

\begin{figure}[ht]
\pic{-0pt}{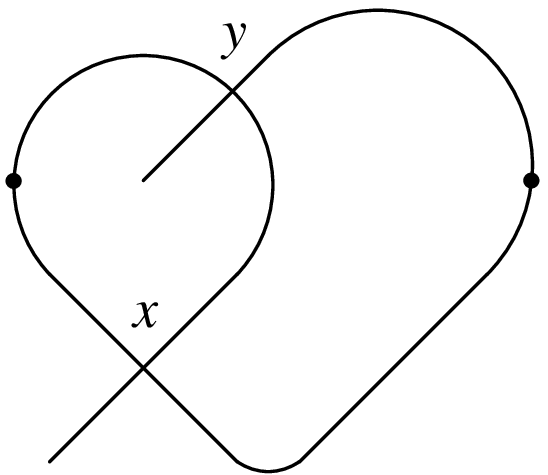}
\caption{}
\end{figure}

There are $11^2$ patterns of signs, represented by a symbol at each double
point. Knots associated with different patterns of signs need not be different. It is still remarkable that we obtain 24 different prime knots in this way. They are listed in Table 2.
The second and third column of Table 2 show the Dowker-Thistlethwaite
numbering and the Rolfsen numbering, respectively. The fourth column tells us
whether the knot is a free divide knot or not. In \cite{GI}, Gibson and Ishikawa
have listed knots of free divides. Up to 10 crossings, their list is complete.
We add the 4-genus in the fifth column. It equals the clasp number and, except
for the knots $9_{46}$, $10_{140}$ and $11n139$, also the Gordian number.

\begin{longtable}{| c | c | c | c | c |}
\caption{Knots associated with a special immersed interval} \\
\hline
\endhead
\hline
$(x,y)$  &  DT numbering  &  Rolfsen numbering  &  free divide  &  $g^*$  \\
\hline 
$(b,c)$  &  $7a4$  &  $7_2$  &  No  &  $1$  \\ \hline
$(b,c_1)$  &  $5a1$  &  $5_2$  &  Yes  &  $1$  \\ \hline
$(b,d)$  &  $7a5$  &  $7_3$  &  Yes  &  $2$  \\ \hline
$(b_1,b_1)$  &  $9n5$  &  $9_{46}$  &  No  &  $0$  \\ \hline
$(b_1,b_3)$  &  $10n29$  &  $10_{140}$  &  No  &  $0$  \\ \hline
$(b_1,c)$  &  $12n121$  &  $-$  &  No  &  $1$  \\ \hline
$(b_1,c_1)$  &  $3a1$  &  $3_1$  &  Yes  &  $1$  \\ \hline
$(b_1,d)$  &  $10n14$  &  $10_{145}$  &  Yes  &  $2$  \\ \hline
$(b_3,b_3)$  &  $11n139$  &  $-$  &  No  &  $0$  \\ \hline
$(b_3,c_3)$  &  $10n4$  &  $10_{133}$  &  No  &  $1$  \\ \hline
$(c,c_3)$  &  $8a2$  &  $8_{15}$  &  No  &  $2$  \\ \hline
$(c,d)$  &  $10n30$  &  $10_{142}$  &  No  &  $3$  \\ \hline
$(c_1,b_1)$  &  $8n2$  &  $8_{21}$  &  No  &  $1$  \\ \hline
$(c_1,b_3)$  &  $9n2$  &  $9_{45}$  &  No  &  $1$  \\ \hline
$(c_1,c_1)$  &  $5a2$  &  $5_1$  &  Yes  &  $2$  \\ \hline
$(c_1,c_3)$  &  $7a3$  &  $7_5$  &  Yes  &  $2$  \\ \hline
$(c_1,d)$  &  $10n31$  &  $10_{161}$  &  Yes  &  $3$  \\ \hline
$(c_3,b_3)$  &  $10n19$  &  $10_{131}$  &  No  &  $1$  \\ \hline
$(c_3,d)$  &  $10n22$  &  $10_{128}$  &  No  &  $3$  \\ \hline
$(d,b_1)$  &  $11n118$  &  $-$  &  No  &  $2$  \\ \hline
$(d,b_3)$  &  $12n407$  &  $-$  &  No  &  $2$  \\ \hline
$(d,c_1)$  &  $7a7$  &  $7_1$  &  Yes  &  $3$  \\ \hline
$(d,c_3)$  &  $10n6$  &  $10_{134}$  &  No  &  $3$  \\ \hline
$(d,d)$  &  $12n591$  &  $-$  &  ?  &  $4$  \\ \hline
\end{longtable}

\medskip
\begin{ack} I am deeply grateful to Daniel Loss and Hans\-peter
Kraft for their great support. Moreover I would like to thank Masaharu
Ishikawa, Stephan Wehrli, Alexander Shumakovitch and Norbert A'Campo for 
motivating me in innumerable discussions.
\end{ack}

\bigskip
\noindent
Mathematisches Institut der Universit\"at Basel

\noindent
Rheinsprung 21

\noindent
CH-4051 Basel

\noindent
Switzerland

\smallskip
\noindent
E-mail address: baader@math-lab.unibas.ch

\end{document}